\documentclass[letterpaper,12pt]{article}
\usepackage{amssymb}
\usepackage{amsmath}
\usepackage{makeidx}
\usepackage{amscd}
\usepackage[all]{xy}

\newtheorem{theo}{Theorem}[section]
\newtheorem{prop}[theo]{Proposition}
\newtheorem{lem}[theo]{Lemma}
\newtheorem{cor}[theo]{Corollary}
\newtheorem{defi}[theo]{Definition}

\def \si {{\sigma}}

\def \Pic {{\rm {Pic\,}}}

\def \Ker {{\rm{Ker\,}}}

\def \A{{\mathbb A}}
\def \P{{\mathbb P}} 
\def \Spec {{\rm{Spec\,}}}
\def \dim {{\rm{dim\,}}}

\def \Hom {{\rm {Hom}}}

\def \Pic {{\rm {Pic}}}
\def \wt {{\rm {wt}}}
\def \GL {{\rm {GL}}}

\def\ov{\overline}

\def \Z {{\mathbb Z}} 
\def \Q {{\mathbb Q}}

\def \WW {{\rm W}}
\def \AA {{\rm A}}
\def \BB {{\rm B}}
\def \DD {{\rm D}}
\def \EE {{\rm E}}
\def \RR {{\rm R}}
\def \CC {{\rm C}} 
\def \FF {{\rm F}} 
\def \GG {{\rm G}}
\def \TT {{\rm T}}

\def \U {{\cal U}}

\def\X{{\cal X}}
\def\Y{{\cal Y}}
\def\G{{\mathbb G}}

\def\T{{\cal T}}
\def\sZ{{\cal Z}} 

\def\lra{\longrightarrow}
 
\def\N{{\cal N}}
\def\H{{\rm H}}
\def\ad{{\rm ad}}
\def\mod{{\rm mod}}

\def\O{{\cal O}}
\def\X{{\cal X}}
\def\Y{{\cal Y}}
 
\def\sZ{{\cal Z }}

\def\O{{\cal O}}

\def\U{{\cal U}}

\def\si{\sigma}

\def\exp{{\rm exp}}
\def\g{{\mathfrak g}}
\def\p{{\mathfrak p}}
\def\n{{\mathfrak n}}
\def\h{{\mathfrak h}}
\def\b{{\mathfrak b}}

\newcommand{\bthe}{\begin{theo}}
\newcommand{\ble}{\begin{lem}}
\newcommand{\bpr}{\begin{prop}}
\newcommand{\bco}{\begin{cor}}
\newcommand{\bde}{\begin{defi}}
\newcommand{\ethe}{\end{theo}}
\newcommand{\ele}{\end{lem}}
\newcommand{\epr}{\end{prop}}
\newcommand{\eco}{\end{cor}}
\newcommand{\ede}{\end{defi}}

\title{Adjoint representation of $\EE_8$ and del Pezzo surfaces of degree $1$}
\author{Vera V. Serganova and Alexei N. Skorobogatov} 
\date{} 
\begin{document}
\baselineskip=15pt
\maketitle
\section*{Introduction}

Let $G$ be the split simple Lie group of type $\EE_8$ with Lie algebra $\g$.
Let $X$ be a split del Pezzo surface of degree 1, and let $\T$ be
a universal torsor over $X$.
In this paper we construct an embedding of $\T$ 
into the $G$-orbit of the highest weight vector of the adjoint
representation of $G$ in $\g$. This orbit is the affine cone $(G/P)_a$
(with the zero removed) over the generalized Grassmannian $G/P\subset\P(\g)$.
Let $H\subset G$ be a split maximal torus, and let $T\subset\GL(\g)$ be the
extension of $H$ by the centre of $\GL(\g)$.
The above embedding is equivariant with respect to the
action of $T$ identified with the N\'eron--Severi torus of $X$. 
Moreover, the $T$-invariant 
hyperplane sections of $\T$ corresponding to the 240 roots of $\EE_8$ are
the inverse images of the 240 lines on $X$. 
This extends the main result of \cite{I} to del Pezzo surfaces
of degree 1.

Generalising the blowing-up construction of \cite[Section 4]{I} 
we prove the following result which may be of independent interest. 
Let $\g=\h\oplus(\bigoplus_{\beta\in \RR}\g_\beta)$ 
be a semisimple Lie algebra with a Cartan subalgebra
$\h$ and a root system $\RR$.
Let $\alpha\in\RR$ be a long simple root, and let $V$ be 
the simple $\g$-module whose highest weight $\omega$
is the fundamental weight dual to $\alpha$.
Define a $\Z$-graded Lie algebra structure on $\g$ by setting 
$\h\subset\g_0$, and $\g_\beta\subset\g_n$ if $\beta-n\alpha$
is a linear combination of simple roots other than $\alpha$.
Then there is a natural $\Z$-grading on $V=\oplus_{n\geq 0}V_n$ 
such that $\g_iV_j\subset V_{j-i}$. The subalgebra $\g_0$ is the direct sum
of the 1-dimensional abelian Lie algebra and a semisimple Lie algebra $\g'$.
Each graded component $V_i$ is a $\g'$-module;
moreover, $V_1$ is the simple $\g'$-module of highest weight $-\alpha$.
Let $G$ (resp. $G'$) be the split simply connected semisimple Lie
group whose Lie algebra is $\g$ (resp. $\g'$), and let $H\subset G$
be the Cartan subgroup with Lie algebra $\h$.
The $G$-orbit of the highest weight vector
in $\P(V)$ is the homogeneous space $G/P$, where $P$
is the maximal parabolic subgroup of $G$ defined by $\alpha$. 
Similarly, $G'/P'\subset \P(V_1)$ is the $G'$-orbit of the highest
weight vector in $\P(V_1)$.
Let $G_{\leq -2}\subset G$ be the unipotent
subgroup whose Lie algebra is the nilpotent Lie subalgebra 
$\g_{\leq -2}\subset\g$. Finally, let $H_\omega$ be the 
1-parameter subgroup of the maximal torus $H$ 
such that the kernel of the natural surjection 
$\hat H=P(\RR)\to\hat H_\omega$ is given by $(x,\omega)=0$.
In Theorem \ref{B1} we construct an open subset of $G/P$
invariant under the semi-direct product
$G_{\leq -2}\rtimes H_\omega$ such that the quotient 
is isomorphic to $\P(V_1)$ blown-up at $G'/P'$.

Although we largely follow the same strategy of proof as in \cite{I}
the generalisation from the cases
$\AA_4$, $\DD_5$, $\EE_6$, $\EE_7$ to the case $\EE_8$ 
is far from straightforward.
The root system $\EE_7$ is obtained from $\EE_8$ by removing
$\alpha=\alpha_8$, the simple root 
corresponding to the last node of the longest leg of the Dynkin diagram.
(Here and elsewhere we use Bourbaki's notation.)
A number of difficulties stem from the fact that the simple Lie algebra 
$\g$ of type $\EE_8$ graded by $\alpha_8$
has five non-zero graded components $\g_n$ and 
not three as was the case for $(\AA_4,\alpha_3)$, $(\DD_5,\alpha_5)$,
$(\EE_6,\alpha_6)$ and $(\EE_7,\alpha_7)$, so in our case $G_{\leq -2}$
is no longer trivial.
The main result of Section \ref{two} is Theorem \ref{B2}
applicable whenever the grading of $\g$ has length 5. 
Let $(G'/P')_a$ be the affine cone over $G'/P'$.
Theorem \ref{B2} says that a natural torsor under the multiplicative
group $\G_m$ over
the blowing-up of a subvariety $Z\subset V_1\setminus\{0\}$ at 
$Z\cap (G'/P')_a$ is isomorphic to a locally closed subset of
$(G/P)_a$ provided there exists
a symmetric bilinear form on $\g_{-1}$ with values in $\g_{-2}$
satisfying certain properties. 
This form allows us to construct a section 
of a quotient morphism by the action of $G_{\leq -2}$.
In Section \ref{three} we zoom in on the cases $\EE_7$ and $\EE_8$ and prove
some technical lemmas about these algebras and related homogeneous spaces.
The preparations for the proof proper start in Section \ref{four}, 
where we construct the required symmetric form, which turns out
to be essentially unique. Its construction 
is made possible by the following fact 
(undoubtedly well known to experts): blowing up a point on a del Pezzo
surface of degree 2 one obtains a del Pezzo surface of degree 1  
only if the point does not belong to the union of exceptional curves
and the branch curve of the anti-canonical double covering (Lemma \ref{PQ}).
The proof of the main result of this paper,
Theorem \ref{main}, is finished in Section \ref{five}. 

\section{The blow-up theorem} \label{one}

Throughout the paper we denote by $k$ a field of characteristic $0$
with an algebraic closure $\ov k$.

\medskip

Let $G$ be a split simply connected semisimple group, with 
a Borel subgroup $B$ defined over $k$, and a split maximal torus 
$H\subset B$, $H\simeq\G_{m,k}^r$. These data define
a root system $\RR$ together with a basis of simple roots $\Delta$. 
Let $\WW$ be the Weyl group of $\RR$. 
If $\alpha\in\RR$, then $\alpha^\vee=\frac{2}{(\alpha,\alpha)}\alpha$
is the corresponding coroot.

Let $\alpha\in\Delta$ be a simple root, and 
$\omega$ be the fundamental weight dual to $\alpha$, that is,
$(\omega,\alpha^\vee)=1$, and $(\omega,\beta^\vee)=0$
if $\beta\in\Delta\setminus\{\alpha\}$.

Let $G\to \GL(V)$ the irreducible representation with the
highest weight $\omega$. Let $P\subset G$, $P\supset B$, 
be the maximal parabolic subgroup such that $G/P$ is the orbit of the 
highest weight vector $v$ in $\P(V)$. The orbit $Gv$ is 
$(G/P)_a\setminus\{0\}$, where $(G/P)_a$ is the
affine cone over $G/P$. 
Let $\tilde G$ be the reductive subgroup of $\GL(V)$
generated by $G$ and the scalar matrices.

Let $\g$, $\h$, $\b$ be the corresponding Lie algebras. 
A simple root $\alpha\in\Delta$ turns 
$\g=\h\oplus(\bigoplus_{\beta\in \RR}\g_\beta)$
into a graded Lie algebra
$\g=\bigoplus_{n\in\Z}\g_n$, where $\h\subset\g_0$ and
$\g_\beta\subset\g_n$ if 
$n$ is the coefficient of $\alpha$
in the decomposition of $\beta$ into an integral 
linear combination of simple roots.
The subalgebra $\p=\g_{\geq 0}$ is the Lie algebra of $P$. 
The subalgebra 
$\g_0$ is reductive, and is the direct sum of the 1-dimensional
centre and the semisimple Lie algebra $\g'=[\g_0,\g_0]$. 
The Dynkin diagram of $\g'$
is obtained from that of $\g$ by removing the node corresponding to $\alpha$.
Let $G'\subset G$ be the semisimple
simply connected group with Lie algebra $\g'$. 

The vector space $V$ is the direct sum $V=\oplus_{n\geq 0} V_n$,
where $V_n$ is spanned 
by the vectors of weight $\tau$ such that $n$ is the coefficient
of $\alpha$ in the decomposition of the root
$\omega-\tau$ into a linear combination of simple roots. 
It is clear that $V$ is a
graded $\g$-module, that is, $\g_iV_j\subset V_{-i+j}$.
We have $V_0=k v$. 

\ble \label{L}
The map $g\mapsto gv$ is an isomorphism 
of $\g'$-modules $\g_{-1}\to V_1$. Moreover,
$V_1$ is an irreducible $\g'$-module with highest weight
$-\alpha$.

The $\g'$-module $V_2$ is isomorphic to $V_2^+\oplus V_2^-$,
where the map $g\mapsto gv$ is an isomorphism of $\g'$-modules
$\g_{-2}\to V_2^-$, and $V_2^+=S^2(V_1)/V(-2\alpha)$,
where $V(-2\alpha)$ is the irreducible 
$\g'$-module with highest weight $-2\alpha$.
\ele
{\em Proof} Let $\U(\g)$ be the universal enveloping algebra of $\g$.
Consider the generalised Verma module $M=\U(\g)\otimes_{\U(\p)}kv$.
By the Poincar\'e--Birkhoff--Witt theorem
the composite map $\U(\g_{\leq -1})\to \U(\g)\to M$ 
is an isomorphism of left $\U(\g_{\leq -1})$-modules, and also of
$\g'$-modules. The grading on 
$\U(\g_{\leq -1})$ induced by the grading on $\g_{\leq -1}$
defines a grading $M=\oplus_{n\geq 0}M_n$.
We have the following decompositions of $\g'$-modules:
$$M_0=kv, \ M_1=\g_{-1}v, \ M_2=\g_{-2}v\oplus S^2(\g_{-1})v.$$
Let $X_{-\alpha}\subset \g_{-\alpha}$ be a non-zero element.
The $\g$-module $V$ is isomorphic to the quotient $M/N$, where the $\g$-submodule
$N$ is generated by $X_{-\alpha}^2v$. 
The standard relations in $\U(\g)$ imply that
$N=\U(\g_{\leq 0})X_{-\alpha}^2v$. The grading on $M$ induces
the grading $N=\oplus_{i\geq 0}N_i$. We have $N_0=N_1=0$, and hence
$V_1=M_1\simeq\g_{-1}$. 
If $\beta\not=\alpha$ is a simple root, then
$\beta-\alpha$ is not a root, thus $X_{-\alpha}$ is 
a highest weight vector of the $\g'$-module $V_1$;
in particular, $V_1$ is an irreducible $\g'$-module
with highest weight $-\alpha$.
The $\g'$-module $N_2$ is generated by $X_{-\alpha}^2v$, thus 
$N_2\simeq V(-2\alpha)$. We obtain
$V_2=M_2/N_2=\g_{-2}\oplus \big(S^2(\g_{-1})/V(-2\alpha)\big)$. 
QED

\medskip

We shall identify the $\g'$-modules $\g_{-1}$ and $V_1$ 
by the isomorphism that sends 
$g$ to $gv$. The exponential map $\exp(x)=\sum_{n\geq 0}\ad(x)^n/n!$
on the nilpotent Lie subalgebra
$\g_{\leq -1}$ is a morphism of affine varieties $\exp:\g_{\leq -1}\to \GL(\g)$
whose image is contained in $G$.
For $x\in \g_{-1}=V_1$ we write 
$$\exp(x)v=v+x+p_2(x)+p_3(x)+\ldots,$$ where $p_n(x)$ is in
$\Hom_{G'}(S^n(V_1),V_n)$. Let $p(x,y)$ be the polarisation of $p(x)=p_2(x)$.
Then $p(x,y)=\frac{1}{2}(xy+yx)v\in V_2^+$
is the symmetric part of $xyv$. The skew-symmetric part of $xyv$ is 
$\frac{1}{2}[x,y]v\in V_2^-$. Note that $p(x)\in V_2^+$.

\medskip

Let $P'\subset G'$ be
the stabiliser of $X_{-\alpha}v$ in $\P(V_1)$. This is a
parabolic subgroup of $G'$, and the affine cone $(G'/P')_a$ over $G'/P'$
is $G'(X_{-\alpha}v)\cup \{0\}$.

We now introduce an important subgroup of $\tilde G$. Define
$D\subset\GL(V)$ as the 1-dimensional torus whose element $g_t$,
$t\in k^*$, acts on $V_i$ as multiplication by $t^{1-i}$.
It is easy to see that $D\subset\tilde G$. Indeed,
let $r$ be the positive rational number such that $r\omega$
is a primitive element of the root lattice $Q(\RR)$.
This lattice is identified with the cocharacter lattice
of $H$. Let $H_\omega\subset H$ be the 1-dimensional subtorus
defined by $r\omega\in Q(\RR)$.
Then $D$ is contained in the 2-dimensional torus generated by
the scalar matrices and $H_\omega$, so that $D\subset\tilde G$.

\ble \label{L1}
$(G'/P')_a=(G/P)_a\cap V_1=p^{-1}(0)$
\ele
{\em Proof} Let us prove the first equality.
The tangent space to $x\in (G/P)_a$ is $kx+\g x$. If
$x\in (G'/P')_a\subset V_1$, then
$$\TT _{x,(G/P)_a} \cap V_1=(kx+ \g x) \cap V_1=kx+\g' x=\TT _{x,(G'/P')_a}.$$
Hence $(G'/P')_a$ is an irreducible component
of $(G/P)_a\cap V_1$. On the other hand, the closed set
$(G/P)_a\cap V_1$ is a union of $G'$-orbits, 
but the closure of any non-zero orbit contains the unique closed orbit
$(G'/P')_a$. Hence $(G'/P')_a=(G/P)_a\cap V_1$. 

If $p(x)=0$, then obviously $p_n(x)=0$ for all $n\geq 2$.
Thus $\exp(x)v=v+x$ is in $(G/P)_a$.
Hence $g_t\,\exp(x)v=tv+x$ is also in $(G/P)_a$ for any 
$t\in k^*$. But $(G/P)_a$ is a closed set, so that
the limit point $x\in V_1$ is contained in it.
By the first equality we see that $x$ is actually
in $(G'/P')_a$. On the other hand, $p(X_{-\alpha}v)=0$,
and since $p$ is $G'$-equivariant, $p$ vanishes
on the orbit $G'(X_{-\alpha}v)$, and hence on $(G'/P')_a$. QED

\medskip

Let $B^-\subset G$ be the opposite Borel subgroup, and 
$N^-\subset G$ its unipotent radical; thus
$B^-=N^-H$. Let $\b^-$ (resp. $\n^-$) be the Lie algebra of $B^-$ (resp.
of $N^-$). Then $N^-=\exp(\n^-)$, and 
$$\n^-=\g_{\leq -1}\oplus (\n^-\cap\g_0)\ \subset\ \g_{\leq 0}.$$
The decreasing family of nilpotent subalgebras $\g_{\leq -n}\subset\n^-$,
$n\geq 1$, defines a decreasing family of unipotent subgroups 
$G_{\leq -n}=\exp(\g_{\leq -n})\subset N^-$.

\medskip

Let $\pi_n:(G/P)_a\to V_n$ be the natural projections. 
Let $\pi_2^+$ (resp. $\pi_2^-$)
be the projection to $V_2^+$ (resp. to $V_2^-$).

\medskip

The Bruhat decomposition represents $G/P$ as a disjoint union
of the Bruhat cells $B^- (k v_\mu)\subset \P(V)$,
where $v_\mu\in V$ is a vector of weight $\mu= w(\omega)$, and
$w$ is a coset representative of $\WW$ modulo the Weyl group of $P$.
Since $V_0=kv$ is the trivial $\g_0$-module,
the big (open) cell is $B^-(kv)=N^-(kv)=G_{\leq -1}(kv)$.
The preimage of the big cell in $(G/P)_a$ 
is a dense open subset of $(G/P)_a$ given by $\pi_0(x)\not=0$.

\ble \label{L2}
If $\mu\in \WW(\omega)$ is a weight of $V_n$, $n\geq 2$, then
$\pi_1(B^- v_\mu)=0$.
\ele
{\em Proof} For any $x\in V_n$ we have 
$B^-x\subset \oplus_{i\geq n}V_i$ since $V$ is a graded $\g$-module. QED

\medskip

Let $G_{\leq -2}\rtimes D\subset \tilde G$ be the semidirect product. 
It is clear that it preserves the fibres of $\pi_1:(G/P)_a\to V_1$.

\bpr \label{p1}
If $x\in V_1$, $x\not\in (G'/P')_a$, then 
$\pi_1^{-1}(x)=(G_{\leq -2}\rtimes D)\exp(x)v$.

If $x\in V_1$, $x\in (G'/P')_a\setminus\{0\}$, then 
$\pi_1^{-1}(x)=G_{\leq -1}x\,\cup\, (G_{\leq -2}\rtimes D)\exp(x)v $.
\epr
{\em Proof} Let $y\in\pi_1^{-1}(x)$, $x\not=0$. 
Then $y$ is contained in $B^-v_\mu$
for some $\mu=w(\omega)$. Since $x\not=0$ we have
$v_\mu=v$ or $v_\mu\in V_1$, by Lemma \ref{L2}. In the first case,
after applying an appropriate element $u\in D$, we ensure that $\pi_0
(uy)=v$ and therefore 
$uy$ is in $G_{\leq -1}v=\exp(\g_{\leq -1}) v$. Since $\pi_1(uy)=\pi_1(y)=x$
we see that $\pi_1^{-1}(x)=D\exp(x+\g_{\leq -2})v=(G_{\leq -2}\rtimes D)\exp(x)v$.
In the second case $y\in V_{\geq 1}$, moreover
$$y\in \exp(\g_{\leq -1})\exp(\n^-\cap \g_0)v_\mu
\subset G_{\leq -1}(G'/P')_a,$$
since $\exp(\n^-\cap \g_0)\subset G'$ and $v_\mu\in (G'/P')_a$.
Now $\pi_1(y)=x$ implies that $x\in (G'/P')_a$ and $y\in G_{\leq -1}x$.
Since $(G'/P')_a$ is a subset of $(G/P)_a$,
we see that $G_{\leq -1}x$ is also a subset of $(G/P)_a$.
This completes the proof. QED

\medskip

It follows that $\pi_1^{-1}(V_1\setminus\{0\})$ is the union of 
$$(G_{\leq -2}\rtimes D)\exp(\g_{-1}\setminus\{0\})v=(G_{\leq -1}\setminus G_{\leq -2})k^*v$$
and $G_{\leq -1}((G'/P')_a\setminus\{0\})$.

{From} now on {\it we assume that $\alpha$ is a {\bf long} 
root of the root system $\RR$}.

\ble \label{L3}
The group $G_{\leq -2}$ acts freely on $\pi_1^{-1}(x)$ for any $x\in V_1\setminus\{0\}$.
\ele
{\em Proof} Recall that $v$ is a vector of highest weight $\omega$, so we
can write $v=v_\omega$.
By Lemma \ref{L2}, $\pi_1^{-1}(V_1\setminus\{0\})$ 
is contained in the union
of $B^-k^*v_\omega$ and $B^-k^*v_\mu$, where $\mu\in\WW\omega$
is a weight of $V_1$, hence it is enough to prove that $G_{\leq -2}$
acts freely on these cells.
If $r_{\alpha}$ is the reflection in the simple root $\alpha$, then
$r_\alpha(\omega)=\omega-\alpha$ is the weight of $X_{-\alpha}v\in V_1$,
thus in the latter case $\mu\in \WW'(\omega-\alpha)$,
where $\WW'$ is the Weyl group of $\g'$. Due to $G'$-invariance it
suffices to check that the stabilisers of $v_{\omega}$ and
$v_{\omega-\alpha}$ in $G_{\leq -2}$ are trivial. Since $G_{\leq -2}$
is unipotent this is equivalent to the triviality of the
stabilisers in the Lie algebra $\g _{\leq -2}$. The
stabiliser of any weight vector $v_{\mu}$ in $\g_{\leq -2}$ is a direct 
sum of root spaces. On the other hand, if $\mu$ is an extremal weight
of $V$ and $\beta$ is a root of $\g$, then either
$\g_{\beta}v_{\mu}=0$ or $\g_{-\beta}v_{\mu}=0$. A simple $sl_2$
argument shows that if $(\mu,\beta)<0$ then
$\g_{-\beta}v_{\mu}=0$ and $\g_{\beta}v_{\mu}\neq 0$. We claim that
$(\mu,\beta)<0$ for $\mu=\omega$ or $\mu=\omega-\alpha$ and 
any root $\beta$ of $\g_{\leq -2}$. Indeed if $\mu=\omega$, then 
$(\omega,\beta)<0$ for any root $\beta$ of $\g_{\leq -1}$.
Now let $\mu=\omega-\alpha$. Then we have
$(\omega-\alpha,\beta)=(r_{\alpha}(\omega),\beta)=(\omega,r_{\alpha}(\beta))$.
Our assumption that $\alpha$ is a long root implies that
$r_{\alpha}(\beta)\in\{\beta-\alpha,\beta,\beta+\alpha\}$, thus $r_{\alpha}(\beta)$
is a root of $\g_{\leq -1}$, so that $(\omega,r_{\alpha}(\beta))<0$.
This implies that $(\omega-\alpha,\beta)<0$
and so completes the proof of the lemma. QED

\bthe\label{B1}
Let $\pi=(\pi_1,\pi_2^+):(G/P)_a\to V_1\oplus V_2^+$ be the natural projection.
Define the open subset $U\subset (G/P)_a$ as the complement to
the union of closed subsets $\pi_1^{-1}(0)$ and $(\pi_2^+)^{-1}(0)$.

{\rm (i)} $G_{\leq -2}$ acts freely on $U$, and the fibres of 
$\pi$ contained in $U$ are orbits of $G_{\leq -2}$.

{\rm (ii)} $G_{\leq -2}\rtimes D$ acts freely on $U$ preserving 
the fibres of 
$U\to V_1\setminus\{0\}\,\times \,\P(V_2^+)$, which are orbits of 
$G_{\leq -2}\rtimes D$.

{\rm (iii)} $G_{\leq -2}\backslash U\to (G_{\leq -2}\rtimes D)
\backslash U$ is a torsor under $D\cong\G_m$.

{\rm (iv)} $(G_{\leq -2}\rtimes D)\backslash U$ is isomorphic to $V_1\setminus\{0\}$
blown up at $(G'/P')_a\setminus\{0\}$. The exceptional divisor is
given by $\pi_0(x)=0$.
\ethe

We write various quotient morphisms in the theorem as a commutative diagram:
\begin{equation}\xymatrix{
U\  \ar[d]^\pi \ar@{^{(}->}[r]&
V_0\times (V_1\setminus\{0\})\times 
(V_2^+\setminus\{0\})\times V_2^-\times V_{\geq 3} \ar[d]^\pi\\
G_{\leq -2}\backslash U \ \ar[d] \ar@{^{(}->}[r]&
(V_1\setminus\{0\})\times (V_2^+\setminus\{0\})\ar[d]\\
(G_{\leq -2}\rtimes D)\backslash U\ \ar@{=}[d] \ar@{^{(}->}[r]&
(V_1\setminus\{0\})\times \P(V_2^+) \ar[d]\\
{\rm Bl}_{(G'/P')_a\setminus\{0\}}(V_1\setminus\{0\}) \ar[r]&V_1\setminus\{0\}}
\label{d1} \end{equation}

\noindent{\em Proof} 
If $t\in k^*$ and $h\in G_{\leq -2}$ are such that
$g_t h\xi=\xi$, then 
$\pi_2^+(\xi)=\pi_2^+(g_t h\xi)=t^{-1}\pi_2^+(\xi)\not=0$,
hence $t=1$. Then $g=1$ by
Lemma \ref{L3}, so that $G_{\leq -2}\rtimes D$ acts freely on $U$.
By Proposition \ref{p1} we can write $U=U_1\cup U_2$, where
$$U_1\subset (G_{\leq -2}\rtimes D)\exp(\g_{-1}\setminus\{0\})v,\ \
\text{and} \ \ 
U_2\subset G_{\leq -1}\big((G'/P')_a\setminus\{0\}\big),$$ since
for $x\in (G'/P')_a$ we have $p(x)=0$ so that
no point in $(G_{\leq -2}\rtimes D)\exp(x)v$ is in $U$.
If $\xi\in U_1$,
then $\pi_2^+(\xi)$ is proportional to $p(\pi_1(\xi))$, thus
Lemma \ref{L1} implies that $\pi_2^+(\xi)\not=0$ is
equivalent to $\pi_1(\xi)\notin (G'/P')_a$, so that
$$U_1=G_{\leq -2}\exp(\g_{-1}\setminus (G'/P')_a)k^*v.$$
The fibres of $\pi$
contained in $U_1$ are orbits of $G_{\leq -2}$, and those of
$U_1\to V_1\setminus\{0\}\,\times \, \P(V_2^+)$ 
are orbits of $G_{\leq -2}\rtimes D$. Moreover,
the morphism 
$$(G_{\leq -2}\rtimes D)\times(V_1\setminus (G'/P')_a)  \ \lra U_1$$
sending $(s,x)$ to $s\,\exp(x)v$, is an isomorphism.
In particular, $\pi_1$ gives rise to a trivial 
$(G_{\leq -2}\rtimes D)$-torsor $U_1\to V_1\setminus (G'/P')_a$. 
Any element of $U_2$ can be written as
$\xi=h\,\exp(y)x$, where $x\in (G'/P')_a\setminus\{0\}$,
$y\in \g_{-1}\simeq V_1$, $h\in G_{\leq -2}$. 
Then $\pi_2^+(\xi)=p(x,y)$, so that
$$U_2=G_{\leq -2}\{\exp(y)x\,|\,x\in(G'/P')_a,\ y\in V_1, \ p(x,y)\not=0\}.$$
Let $x$ be a non-zero point of $(G'/P')_a$. 
Let us observe that $p(x,y)=0$ for $y\in V_1$ if and only if 
$y$ is in the tangent space $\TT _{x,(G'/P')_a}$, 
since $p(x)=0$ gives a system
of quadratic equations defining $(G'/P')_a$, by Lemma \ref{L1}.
Thus the zero set of $p(x,y)$ in 
$\big((G'/P')_a\setminus\{0\}\big)\times V_1$
is the tangent bundle of $(G'/P')_a\setminus\{0\}$.
Moreover, for such pairs $(x,y)$ we have $\exp(y)x=x$. 
For this we need to show that $yx=0$, and this follows from
$[x,y]=0$ by the remarks after Lemma \ref{L}, so we
only need to prove that $x$ and $y$ commute.
Recall that $\TT _{x,(G'/P')_a}$ is $kx+\g'x\subset V_1$. 
By the $G'$-invariance we can assume without loss of generality
that $x=X_{-\alpha}v$, so that we must show that $[X_{-\alpha},[X_{-\alpha},\g']]=0$.
For this it is enough to prove that $[X_{-\alpha},[X_{-\alpha},X_\beta]]=0$
for any root $\beta$ of $\g'$. But $\beta\not=\alpha$, and it is well known
that $\beta-2\alpha$ is never a root for any long root $\alpha\not=\beta$.
This finishes the proof that $\exp(y)x=x$.

Let us show that the fibres of the restriction of $\pi$ to $U_2$
are orbits of $G_{\leq -2}$. If $\exp(y)x$ and $\exp(y')x'$
have the same image under $\pi$,
then $x'=x$ and $p(x,y)=p(x,y')$, so that $y'-y\in \TT _{x,(G'/P')_a}$.
As we have seen, this implies $\exp(y'-y)x=x$. Since
$$\exp(y')=h\,\exp(y)\exp(y'-y)$$ 
for some $h\in G_{\leq -2}$,
we are done. It follows that the fibres of 
$U_2\to V_1\setminus\{0\}\,\times \, \P(V_2^+)$ 
are orbits of $G_{\leq -2}\rtimes D$,
which completes the proof of (i) and (ii). Part (iii) is now obvious.

Let $\N$ be the normal bundle to $(G'/P')_a\setminus\{0\}$ in $V_1$, 
that is, the cokernel of the injective map of vector bundles 
$\TT _{(G'/P')_a}\to V_1$.
The map $(x,y)\mapsto (x,p(x,y))$ identifies 
$\N$ without its zero section with 
$$G_{\leq -2}\backslash U_2\,\subset\,
\big((G'/P')_a\setminus\{0\}\big)\,\times\, (V_2^+\setminus\{0\}),$$
thus $(G_{\leq -2}\rtimes D)\backslash U_2=\P(\N)$. 
Finally, $\pi_1:(G_{\leq -2}\rtimes D)\backslash U\to V_1\setminus\{0\}$
is a morphism of smooth varieties which is an isomorphism 
away from $(G'/P')_a$, whereas $\pi_1^{-1}((G'/P')_a\setminus\{0\})$
is isomorphic to the projectivisation of the normal bundle to 
$(G'/P')_a\setminus\{0\}$ in $V_1\setminus\{0\}$. 
It is known and not very hard to
prove that this implies the first statement of (iv). 
But $U_2$ is the closed subset of $U$ given by
$\pi_0(x)=0$. This finishes the proof.
QED

\section {The case of grading of length 5} \label{two}

Let us now assume that the grading of $\g$ 
defined by a simple root $\alpha$ has length 5, i.e.,
$\g_n=0$ exactly when $|n|>2$. An inspection of tables in \cite{B}
shows that this is the full list of such pairs $(\RR,\alpha)$:
$$(\BB_n,\alpha_i),\, i\not=1; \ 
(\CC_n,\alpha_i),\, i\not=n; \ 
(\DD_n, \alpha_i),\, i\notin\{1,n-1,n\}; \
(\EE_6,\alpha_i),\, i=2,\,3,\,5;$$
$$(\EE_7,\alpha_i),\, i=1,\,2,\,6; \
(\EE_8,\alpha_i),\, i=1,\,8; \ 
(\FF_4,\alpha_i),\, i=1,\,4;\  (\GG_2,\alpha_2).$$ 
Recall that our
enumeration of roots follows the conventions of \cite{B}.

We keep the notation of the previous section,
in particular $V$ is the simple $\g$-module with highest weight $\omega$,
the fundamental weight dual to $\alpha$. We identify $V_1$ with
$\g_{-1}$, and $V_2^-$ with $\g_{-2}$.

\bthe\label{B2}
Assume that the grading of $\g$ defined by a simple root $\alpha$
has five non-zero terms.
Let $Z$ be a smooth closed subset of $\g_{-1}\setminus\{0\}$ 
such that $Z_0:=Z\cap(G'/P')_a$ is also smooth.
Assume that there exists a linear map $s:S^2(\g_{-1})\to \g_{-2}$ 
such that $s(x)=0$ and $[a,x]=4 s(x,a)$ for any 
$x\in Z_0$ and $a\in \TT _{x,Z}$.

Define $\tilde{Z}=D\{\exp (x+s(x))v | x\in Z\}\cap U$, and 
let $\sZ$ be the Zariski closure of
$\tilde{Z}$ in $\pi_1^{-1}(Z)\cap U$. Then

{\rm (i)} $\pi_1:\sZ\to Z$ is surjective.

{\rm (ii)} $\sZ$ is $D$-invariant, and $D$ acts freely on $\sZ$.

{\rm (iii)} The quotient $D\backslash \sZ$ is isomorphic to $Z$
blown up at $Z_0$. The exceptional divisor is
given by $\pi_0(x)=0$.
\ethe
This theorem states that the above sets are related as follows:
$$\xymatrix{
\tilde Z\, \ar@{^{(}->}[r]  \ar@<1ex>[d]^{\pi_1}&
\sZ \ar@<1ex>[d]_{\pi_1} \ar[r]& D\backslash \sZ \ar@{=}[d] \\
Z\setminus Z_0 \ar@<1ex>[u]^{\exp(x+s(x))v}\ar@{^{(}->}[r]&Z&
{\rm Bl}_{Z_0}(Z)\ar[l]}$$
where the downward arrows $\pi_1$ are surjective.
\medskip

\noindent{\em Proof} {\rm (i)} If $x\in Z\setminus Z_0$, then 
$\exp(x+s(x))v\in U$ because $x\not=0$ and $p(x)\not=0$
by Lemma \ref{L1}. Thus 
$\exp(x+s(x))v\in \tilde Z\subset \sZ$. Since
$x=\pi_1(\exp(x+s(x))v)$, we see that $x\in \pi_1(\sZ)$.

Let $k[[t]]$ be the $k$-algebra of formal power series.
Now let $x\in Z_0$ and $a\in \TT _{x,Z}$, and let 
$$\phi(t)=x+at+O(t^2)\in Z(k[[t]])$$
be a $k[[t]]$-point of $Z$. Let us prove that 
$$y=\lim_{t\to 0}\, g_t\, \exp \big(\phi(t)+s(\phi(t))\big)v$$
is a well defined point of $\sZ$.
Using the identity
$g\, \exp (h) g^{-1}=\exp({\rm Ad}_g h)$
and the fact that $g_t\,(v)=tv$ we obtain
$$y=\lim_{t\to 0}\, \exp \big(g_t\,(\phi(t)+s(\phi(t)))g_t^{-1}\big) tv.$$
Since ${\rm Ad}_{g_t}z=t^i z$
for any $z\in \g_i$, and
$$\phi(t)+s(\phi(t))=x+at+O(t^2)+2s(x,a)t+O_2(t^2),$$
where $O_2(t^2)\in \g_{-2}$, we obtain
$$y=\lim_{t\to 0}\,\exp \big(xt^{-1}+a+O(t)+2s(x,a)t^{-1}+O_2(1)\big)tv.$$
Since $[\g_{-2},\g_{\leq -1}]=0$,
by the Campbell--Hausdorff formula for any $b,\,c\in\g_{\leq -1}$ we have
$$\exp (b) \exp(c)=\exp \big(b+c+\frac{[b,c]}{2}\big).$$
Since $O_2(1)\in \g_{-2}$ and $4s(x,a)=[a,x]$ we have
$$y=\lim_{t\to 0}\,\exp(O(1))\exp(a) \exp(xt^{-1})tv=
\lim_{t\to 0}\,\big(\exp(a)x v+O(1)xtv\big)=\exp(a)x v,$$
where we used that $[x,[x,v]]=2p(x)=0$ which holds because $x$ is in
$Z_0$. Thus, $y$ is well defined and, moreover,
\begin{equation}\label{bl}
y=\exp(a)x v=xv+p(x,a)v \ \mod\, V_2^-\oplus V_{>2}.
\end{equation} 
In particular $\pi_1(y)=x$. Hence $\pi_1:\sZ\to Z$ is surjective.

{\rm (ii)} follows from the $D$-invariance of $\tilde{Z}$ and Theorem 
\ref{B1}(ii).

{\rm (iii)} Let $\Y=\pi_1^{-1}(Z)\cap U$ and $\X=\pi(\Y)$.
It is clear that $\Y$ is a closed subset of $\X\times_k\A^n_k$,
where $\A^n_k=V_0\oplus V_2^-\oplus V_{\geq 3}$. 
By construction $\sZ$ is closed in $\Y$, and hence is closed in 
$\X\times_k\A^n_k$:
$$\xymatrix{
\sZ\, \ar@{^{(}->}[r]  \ar[dr]&
\Y\, \ar@{^{(}->}[r] \ar[d]_{\pi} &\X\times_k\A^n_k \ar[dl] \\
&\X&}$$
We shall prove that $\pi$ induces an isomorphism 
$\sZ\tilde\lra \X$. 
By the functoriality of blowing up and Theorem \ref{B1} (iv),
$(G_{-2}\rtimes D)\backslash\Y\simeq D\backslash\X$ is
isomorphic to ${\rm Bl}_{Z_0}(Z)$, so this is enough to complete
the proof of (iii).

Write $\X_0=\pi(\pi_1^{-1}(Z_0))$. We have the 
following useful descriptions of $\X_0$ and its complement in $\X$:
$$\X_0=\{(x,p(x,a))\in U | x\in Z_0,a\in T_{x,Z}\},$$
$$\X\setminus \X_0=\{(x,t p(x)) | x\in Z\setminus Z_0, t\in k^*\}.$$
The image $\pi(\sZ)$ contains $\X\setminus \X_0$ by the argument
from the beginning of the proof of (i), and it contains $\X_0$ by
formula (\ref{bl}). Thus $\pi(\sZ)=\X$.

Let us show that $\pi$ induces an isomorphism 
$\pi^{-1}(\X\setminus \X_0)\cap\sZ\tilde\lra \X\setminus \X_0$.
If $z=(z_0,z_1,z_2^+,z_2^-,...)$ is a $\ov k$-point of $\sZ$,
then we have
\begin{equation}
z_0z_2^+=p(z_1),\quad z_0z_2^-=s(z_1)v, \label{j}
\end{equation}
because these equations are satisfied on the open subset 
$\tilde Z\subset\sZ$ which is given by $z_0\not=0$.
If $y=(y_1,y_2^+)$ is a $\ov k$-point of $\X\setminus\X_0$, 
then $p(y_1)\neq 0$ and $y_2^+\not=0$. Suppose that $\pi(z)=y$, then
$z_1=y_1$ and $z_2^+=y_2^+$.
That implies that $z=g_t\,\exp(y_1+s(y_1))v$, where $t$ is such that
$y_2^+ t=p(y_1)$,
is a unique point of $\sZ$ above $y$. This defines
a section of $\pi:\sZ\to \X$ over $\X\setminus \X_0$. 
Applying Lemma \ref{affine} below with $A=\sZ$ and $B=\X$
we conclude that $\X$ is isomorphic to $\sZ$ . 
The second statement of (iii) is immediate from the first
equation of (\ref{j}). QED

\ble\label{affine}
Let $B$ be a normal 
geometrically integral variety over a field $k$, and
let $\varphi$ be the projection $B\times_k\A^n_k\to B$.
Let $A\subset B\times_k\A^n_k$ be a closed irreducible subscheme
such that $\varphi(A)=B$. If $\varphi$ induces an isomorphism
of fields of functions $k(B)\tilde\lra k(A)$, then
$\varphi: A\to B$ is an isomorphism.
\ele
{\em Proof} Let us denote the field $k(A)=k(B)$ by $K$. 
Let $\Omega\subset B$ be the largest open subset such that 
$\varphi$ induces an isomorphism 
$\varphi^{-1}(\Omega)\cap A\tilde\lra\Omega$.
Let us show that $B\setminus\Omega$ has codimension at least 2,
i.e., $\Omega$ contains all the points of $B$ of codimension 1.
Let $D\subset B$ be an irreducible divisor, and let $\O_D\subset K$
be its local ring. Since $B$ is normal, $\O_D$ is a discrete valuation
ring with valuation $val:\O_D^*\to\Z$. Write $\Spec(\O_D)\times_B A=\Spec(R)$,
where $R$ is a subring of $K$ that contains $\O_D$. 
If $val(x)< 0$ for some $x\in R\setminus 0$, then
$R=K$ and the closed fibre of $\Spec(R)\to \Spec(\O_D)$ is empty. 
Since $\Spec(R)\to \Spec(\O_D)$ is surjective we conclude that
$val(x)\geq 0$ for all $x\in R\setminus 0$, hence
$R=\O_D$. Therefore, the codimension of 
$B\setminus\Omega$ is at least 2. The composition of 
$\varphi^{-1}:\Omega\to A$ with any coordinate projection 
$A\subset B\times_k\A^n_k\to \A^1_k$ is a rational function on $B$
which is regular away from a closed subset of codimension 2, and hence
is regular everywhere on $B$. Since $A$ is irreducible we have 
$\varphi^{-1}(B)=A$, so that $\varphi$ is indeed an isomorphism. QED

\medskip

We thank J-L. Colliot-Th\'el\`ene for pointing out this simple proof.

\section{The case when the adjoint representation is fundamental} \label{three}

Consider the case when the adjoint representation of $\g$
is a fundamental representation, i.e., when the maximal root of $\RR$ 
is the fundamental weight $\omega$ dual to some simple root $\alpha$. 
This happens precisely in the following cases:
$$(\BB_n,\alpha_2),\ n\geq 3; \ 
(\DD_n, \alpha_2),\ n\geq 4; \
(\EE_6,\alpha_2),\ (\EE_7,\alpha_1),\ (\EE_8,\alpha_8),\ 
(\FF_4,\alpha_1),\  (\GG_2,\alpha_2).$$ 
The tables in \cite{B} show that
the coefficient of $\alpha$ in the decomposition of 
the root $\omega$ into a linear combination of simple roots is 2.
Thus the $\mathbb Z$-grading $\g=\oplus \g_n$ defined by $\alpha$
has exactly five non-zero terms $\g_n$, $|n|\leq 2$.
The following properties are easy to check:
\begin{itemize}
\item $\g_0=\g'\oplus  kz$ is the direct sum of 
Lie algebras, where $z\in\h$, $z\not=0$, spans the centre of $\g_0$, 
and $\g'$ is semisimple;
\item the $\g'$-modules $\g_{\pm 1}$ are isomorphic symplectic 
irreducible $\g'$-modules such that all weights have multiplicity 1;
\item $\g_{\pm 2}$ are trivial $\g'$-modules, $\dim \g_{\pm 2}=1$;
\item $z$ is a grading element of $\g$, i.e.
$[z,g]=ng$ for any $g\in \g_n$.
\end{itemize}
We can choose generators $v\in \g_2$ and $w\in \g_{-2}$ so that 
$z=[v,w]$. Using that $\g_n=0$ for $|n|\geq 3$ one checks  
that $[[y,w],v]=y$ for any $y\in \g_1$,
and $[[x,v],w]=x$ for any $x\in\g_{-1}$.
We identify the $\g_0$-modules $\g_{-1}$ and $\g_1$ via the isomorphism
that sends $x$ to $[x,v]$; its inverse sends $y$ to $[y,w]$. 
Define a symplectic form on $\g_{-1}$ by
$$\langle a,b\rangle w=[a,b],$$
where $a,\,b\in \g_{-1}$. It is easy to check that this form is non-degenerate.

\ble \label{sym}
For any $y\in (G'/P')_a\subset\g_{-1}$ and 
$a\in \TT _{y,(G'/P')_a}$ we have $\langle y,a\rangle=0$.
\ele
{\em Proof} Recall that $\g_{2}=\g_\omega$, where 
$\omega$, the fundamental weight dual to $\alpha$, 
is the highest weight of the adjoint representation of $\g$.
Recall also that the highest weight of the $\g'$-module $\g_{-1}$
is $\omega-\alpha$. Since 
the symplectic form $\langle a,b\rangle$ is $G'$-invariant,
it is enough to prove the statement when $y\in \g_{-1}$ is an eigenvector
of $H$ of weight $\omega-\alpha$. 
Since $\TT _{y,(G'/P')_a}=ky+[\g',y]$ we must prove that
the vector space $[\g',y]$ has zero intersection with $\g_{\alpha}$.
This follows from the fact
that $\omega-2\alpha$ is not a root. QED

\medskip

Define the invariant tensors
$p\in \Hom_{\g'}(S^2(\g_{-1}),\g_0)$, $q\in \Hom_{\g'}(S^3(\g_{-1}),\g_{-1})$,
$r\in\Hom_{\g'}(S^4(\g_{-1}),k)$ as follows:
$$p(x)=\frac{1}{2}\ad_x^2(v), \quad q(x)=\frac {1}{6} \ad_x^3(v), \quad
r(x)w=\frac {1}{24}\ad_x^4(v).$$
Then for any $x\in \g_{-1}$ we can write $\exp (x)v$
as the sum of graded components
\begin{equation}
\exp(x)v=
v+[x,v]+p(x)+q(x)+r(x)w. \label{e1}
\end{equation}
We denote the polarisations of these forms by the same letters,
for example
$$
r(a,b,c,d)=\frac{1}{576}\sum_{\pi\in S_4}\ad_{\pi(a)}\ad_{\pi(b)}\ad_{\pi(c)}\ad_{\pi(d)}(v).
$$

\ble \label{1}
For any $x\in \g_{-1}$ we have $\ad_x^2(v)=2p(x)\in\g'$. 
\ele
{\it Proof} 
For any $x\in \g_{-1}$ we have $[[[x,v],x],v]=0$ from the Jacobi identity,
hence $[x,[x,v]]\in\g'$ which proves our formula.  QED

\medskip

The intersection $\h'=\g'\cap\h$ is a Cartan subalgebra in $\g'$.
Since $\g_{-1}$ is a minuscule $\g'$-module,
all the weights of $\g_{-1}$ have multiplicity 1 with respect to $\h'$. 
Let $\Lambda\subset (\h')^*$ be the set of weights of $\g_{-1}$.
Let $X_\mu\in \g_{-1}$ be a non-zero vector of weight 
$\mu\in \Lambda$. Then any $x\in \g_{-1}$
is uniquely written as $x=\sum x^\mu X_\mu$, where $x^\mu$
is a homogeneous coordinate of weight $\mu$. 
Set $c_\mu=\langle X_\mu,X_{-\mu}\rangle$.
Then clearly $c_{-\mu}=-c_\mu$. These numbers are non-zero
since the symplectic form $\langle x,y\rangle$
is non-degenerate. We can write
$$r(x)=\sum_{\mu_1+\mu_2+\mu_3+\mu_4=0}
r_{\mu_1,\mu_2,\mu_3,\mu_4}x^{\mu_1}x^{\mu_2}x^{\mu_3}x^{\mu_4},$$
where the monomials correspond to all sets of four (not necessarily distinct)
elements of $\Lambda$ with zero sum.

Write $q(x)=\sum q^\mu(x) X_\mu$. 
We have $[X_\mu,x]=c_\mu x^{-\mu}w$ and 
$[X_\mu,q(x)]=c_\mu q^{-\mu}(x)w$.

\ble \label{derive}
We have the following formulae:
\begin{equation}
\frac{\partial r(x)}{\partial x^\mu}=c_\mu q^{-\mu}(x)=
\langle X_\mu,q(x)\rangle, \label{r}
\end{equation}
\begin{equation}
\frac{\partial q(x)}{\partial x^\beta}=
[X_\beta,p(x)]+\frac{1}{2}c_\beta x^{-\beta}x=
[X_\beta,p(x)]+\frac{1}{2}\langle X_\beta,x\rangle x.
\label{q}
\end{equation}
\ele
{\em Proof} 
The left hand side of (\ref{r}) multiplied by $24$ is
$$[X_\mu,[x,[x,[x,v]]]]+[x,[X_\mu,[x,[x,v]]]]+[x,[x,[X_\mu,[x,v]]]]+
[x,[x,[x,[X_\mu,v]]]].$$
Here the first term equals $6[X_\mu,q(x)]$. The second term is
$$6[X_\mu,q(x)]+[[x,X_\mu],[x,[x,v]]]=6[X_\mu,q(x)],$$
since $p(x)\in \g'$ by Lemma \ref{1}, and $\g'$ is the stabiliser of $w$.
The third term equals $[x,[X_\mu,[x,[x,v]]]]+[x,[[x,X_\mu],[x,v]]]$,
but since $[w,[x,v]]=-x$, it is the same as the second term. Finally,
the fourth term is 
$[x,[x,[X_\mu,[x,v]]]]+[x,[x,[[x,X_\mu],v]]]$. Using $[w,v]=-z$
and $[z,x]=-x$ we conclude that it is the same as the third term,
thus completing the proof of (\ref{r}).

The left hand side of (\ref{q}) multiplied by $6$ is the following expression
$$[X_\beta,[x,[x,v]]]+[x,[X_\beta,[x,v]]]+[x,[x,[X_\beta,v]]].$$
The first term equals $2[X_\beta,p(x)]$.
The second term is the sum of the first term and
$[[x,X_\beta],[x,v]]=c_\beta x^{-\beta}x$. The last term is the sum of
the second term and 
$[x,[[x,X_\beta],v]]=-c_\beta x^{-\beta}[x,[w,v]]=c_\beta x^{-\beta}x$.
QED

\medskip

For the sake of completeness we list here some formulae that follow from
(\ref{r}) and (\ref{q}), but which will not be used in the rest of the paper:
$$
\frac{\partial^2 r(x)}{\partial x^\mu\partial x^\beta}w=
[X_\mu,[X_\beta,p(x)]]+
\frac{1}{2}\langle X_\mu,x\rangle\langle X_\beta,x\rangle w; $$
for any $x,\,a\in \g_{-1}$ we have 
$$
q(x,x,a)=\frac{1}{3}[a,p(x)]+\frac{1}{6}\langle a,x\rangle x,
\quad\quad r(x,x,x,a)=\frac{1}{4}\langle a,q(x)\rangle. $$

\medskip

{From} now on {\it we only consider the case when $\g$ is 
the simple Lie algebra of type $\EE_8$.}
We have $\dim(\g)=248$, $\g'$ is the simple Lie algebra of type 
$\EE_7$, $\dim \g'=133$, and $\dim \g_{\pm 1} =56$.
It is known that $\g$ is the algebra of endomorphisms of $\g_{-1}$
preserving the invariant quartic form $r(x)$ and the symplectic 
form $\langle x,y\rangle$.
It is also known that all coefficients 
$r_{\mu_1,\mu_2,\mu_3,\mu_4}$ are non-zero (and can be made 
$\pm 1$ with an appropriate choice of basis vectors $X_\mu$,
see \cite{Lurie}, Thm. 6.1.2).
We deduce from (\ref{r}) that for any $\mu\in \Lambda$ the cubic form
$q^\mu(x)$ is a linear combination of the monomials
$x^{\mu_1}x^{\mu_2}x^{\mu_3}$ such that 
$\mu_1+\mu_2+\mu_3=\mu$ with non-zero coefficients.
In particular, $q^\mu(x)$ is not divisible by $x^\mu$,
for any $\mu\in \Lambda$.

The following technical lemma will be used later in the construction.

\ble \label{aux1} 
Let $\rho(x)\in S^3(\g_{-1}^*)$ be a non-zero homogeneous 
cubic form of weight $\beta\in \Lambda$, 
$$\rho(x)=\sum_{\mu_1+\mu_2+\mu_3=\beta}
c_{\mu_1,\mu_2,\mu_3}x^{\mu_1}x^{\mu_2}x^{\mu_3}.$$ 
If $c_{\mu_1,\mu_2,\mu_3}=0$ whenever
$\mu_i=\beta$, then $\rho(x)$ is not identically zero on $(G'/P')_a$.
\ele
{\em Proof} Let $\g'_0\subset \g'$ denote the stabiliser of the hyperplane
of $\g_{-1}$ given by
$x^{\beta}=0$. Then $\g'_0\simeq k\oplus \g''$ is a direct sum of Lie algebras,
where $\g''$ is the simple Lie algebra of type ${\rm E}_6$. 
The $\g''$-module $\g_{-1}$ is the direct sum of irreducible submodules
\begin{equation}\label{dec}
\g_{-1}=W_{-3}\oplus W_{-1}\oplus W_1\oplus W_3,
\end{equation}
where $W_3$ and $W_{-3}$ are trivial $\g''$-modules of dimension $1$,
$X_\beta\in W_3$, and $W_1$ and $W_{-1}$
are dual $\g''$-modules of dimension 27. Moreover, there exists an
element $h$ in the centre of $\g'_0$ such that $[h,v]=i v$ for any
$v\in W_i$. The polynomial $\rho(x)$ must have weight 3 
with respect to $h$, so
$$\rho(x)\,\in\, (S^2(W_3^*)\otimes W^*_{-3})\, \oplus \, 
(W^*_3\otimes W^*_1\otimes W^*_{-1})\, \oplus\, S^3(W^*_1).$$ 
Let $\phi:\g_{-1}\to W_1$ be the natural projection.
Since $W_3^*$ is spanned by $x^\beta$
our assumption on $\rho(x)$ implies that $\rho(x)\in S^3(W_1^*)$, so that
$\rho(x)=\rho(\phi(x))$ for any $x$ in $\g_{-1}$.
Since $\phi((G'/P')_a)=W_1$, the non-zero form
$\rho(x)$ cannot vanish everywhere on $(G'/P')_a$. QED

\section{Curves on del Pezzo surfaces of degree $2$} \label{four}

For a curve $C$ on a smooth surface $X$ we write $[C]$ for the class
of $C$ in the Picard group $\Pic X$. We denote the intersection index of
divisors $D_1$ and $D_2$ on $X$ by $(D_1.D_2)$.

\ble \label{PQ} Let $M_1,\ldots, M_8$ be points in $\P^2_k$
such that the blow-up of $\P^2_k$ in $M_1,\ldots, M_8$
is a del Pezzo surface $X$ of degree $1$. Let $X'$
be the del Pezzo surface of degree $2$ obtained by blowing-up
$\P^2_k$ at $M_1,\ldots, M_7$, and let $M$ be the point corresponding 
to $M_8$ in $X'$. Let $B\subset X'$ be the branch curve
of the anti-canonical double covering $\kappa:X'\to \P^2_k$.
Then $M\not\in B$. The induced map of cotangent spaces
$$\kappa^*:\TT^*_{\kappa(M),\P^2}\ \lra \ \TT^*_{M,X'}$$
is an isomorphism.
\ele
{\em Proof} Let $\si:X\to X'$ be the morphism inverse 
to the blowing-up of $M$ in $X'$, and let $E=\si^{-1}(M)$
be the exceptional divisor. 
Since $X$ is a del Pezzo surface of degree $1$ it is 
clear that $M$ does not belong to the exceptional curves of $X'$.
It is well known that $B\subset\P^2_k$ is a smooth quartic curve,
and that the union of exceptional curves of $X'$ is
the inverse image of the union of bitangents to $B\subset\P^2_k$, 
see \cite[Ch. 4]{M}. Thus if $M\in B$,
then the tangent line $L$ to $B$ at $\kappa(M)$ is not a bitangent.
Hence $\kappa^{-1}(L)$ is a rational curve with one node and no other singular
points. Let $C$ be the strict transform of $\kappa^{-1}(L)$ in $X$, that is,
the closure of $\kappa^{-1}(L)\setminus\{M\}$ in $X$. The multiplicity
of $M$ in $\kappa^{-1}(L)$ is $2$, hence the intersection index $(C.E)=2$.
For the same reason we have the following relation in $\Pic X$:
$$[C]+2[E]=\si^*(\kappa^{-1}(L))=\si^*(-K_{X'})=-K_X+[E].$$
Hence $[C]=-K_X-[E]$ and so $(C.-K_X)=0$, which contradicts the ampleness
of $-K_X$. QED

\medskip



Let $T'\subset \GL(\g)$ be the torus generated by 
the maximal torus $H'\subset G'$ and the 1-dimensional
torus $\G_m$ whose element $t\in k^*$ acts on $\g_n$
as multiplication by $t^{n+2}$. (Note that $H'\cap \G_m=\{\pm 1\}$.)
We denote by $\chi_0$ the character of $T'$ by which $T'$ acts on the 
1-dimensional centre of $\g_0$. This gives natural exact sequences
$$0\lra Q(\EE_7)\lra \hat T'\lra\widehat{\G}_m=\Z\lra 0$$
and
$$0\lra \Z\chi_0\lra \hat T'\lra \hat H'=P(\EE_7)\lra 0.$$

For $\chi\in \hat T'$ let
$S^n_\chi(\g_{-1})$ be the weight $\chi$ eigenspace of $S^n(\g_{-1})$,
and let $S^n_\chi(\g_{-1}^*)$ be the dual space of 
homogeneous forms. In other words, we have 
$\phi(x)\in S^n_\chi(\g_{-1}^*)$ if and only if
$\phi(tx)=\chi(t)^{-1}\phi(x)$.

Define $(G'/P')_a^{\rm sf}$ as the open subset of $(G'/P')_a$
consisting of stable points with respect to $H'$
(which means that the $H'$-orbits are closed with finite stabilisers), 
with the additional
condition that the stabilisers in $T'$ are trivial. By geometric
invariant theory \cite{GIT} the quotient
$Y'=T'\backslash (G'/P')_a^{\rm sf}$ exists as a quasi-projective
variety. By \cite[Thm. 2.7]{I} the canonical morphism
$f':(G'/P')_a^{\rm sf}\to Y'$ is a universal torsor. 
By \cite[Thm. 6.1]{I} there is an embedding 
$X'\hookrightarrow Y'$ such that the images of the weight hyperplane
sections cut the exceptional curves on $X'$. Moreover,
the restriction of $f'$ to $X'$ is a universal torsor $\T'\to X'$,
and so defines an isomorphism $\hat T'\to\Pic X'$. 
It follows that the natural restriction map
$\Pic Y'\to\Pic X'$ is an isomorphism.
The type of the universal torsor $f':\T'\to X'$ up to sign is an isomorphism
$\tau:\hat T'\tilde\lra \Pic X'$ described on page 397 of \cite{I}.
We reproduce this description here for the convenience of the reader.
Let $\chi\in \hat T'$, and let $\phi(x)\in S^n_\chi(\g_{-1}^*)$.
Let $Z_\phi\subset \T'$ be the closed $T'$-invariant subset given
by $\phi(x)=0$, and let $C_\phi=X'\cap f'(Z_\phi)$.
If $\phi(x)$ does not vanish identically
on $\T'$, then $C_\phi$ is a curve on $X'$ 
whose class in $\Pic X'$ equals $\tau(\chi)$.
Following a convention of \cite{I} we identify $\hat T'$ with $\Pic X'$
via isomorphism $-\tau$. Then
by formula (14) of \cite{I}
the intersection index $(C_\phi.-K_{X'})$, also called the degree of $C_\phi$,
equals $n$. Moreover, by formula (15) of \cite{I} we have 
\begin{equation}
\H^0(X',\O(C_\phi))=k[\T']_{-\chi}=
S^2_{\chi}(\g^*_{-1})/\big(I(\T')\cap S^2_{\chi}(\g^*_{-1})\big).
\label{form}
\end{equation}
Here are some important examples of curves of low degree on $X'$.
If $n=1$ and $\mu\in \hat T'$ is a weight of $\g_{-1}$, 
we denote by $\ell_\mu$ the exceptional curve in $X'$
cut by the image of the hyperplane section given by $x^\mu=0$.
It is clear that $[\ell_\mu]=\mu$.
We note that $\mu$ is a weight of $\g_{-1}$ if and only
if $\chi_0-\mu$ is a weight of $\g_{-1}$.
According to formula (12) of \cite{I} the intersection index
of $\ell_\mu$ and $\ell_\nu$ can be written as
\begin{equation}
(\ell_\mu.\ell_\nu)=\frac{1}{2}-(\mu,\nu), \label{int}
\end{equation}
where the last pairing is the standard bilinear form on
$Q(\EE_7)\otimes\Q$ applied to the restrictions of $\mu$ and $\nu$
to $\hat H'=P(\EE_7)$.

For $n=2$ we have $S^2(\g_{-1})=V_2^+\oplus V(-2\alpha)=\g'\oplus V(-2\alpha)$
(cf. Lemma \ref{L}).
If $\phi(x)\in S^2_{\chi_0}(\g_{-1}^*)$, then $C_\phi$
is an anti-canonical curve, i.e. $[C_\phi]=\chi_0=-K_{X'}$. 
Indeed, this is the only effective divisor class with self-intersection $2$
which is orthogonal to $Q(\EE_7)\subset \Pic X'$.

Now let $\xi\in \hat T'$ be a weight of $\g'$, $\xi\not=\chi_0$. Then
it can be checked using (\ref{int}) that 
$\xi=\mu+\nu$, where $\mu,\,\nu\in \hat T'$
are weights of $\g_{-1}$ such that the intersection index
$(\ell_\mu.\ell_\nu)=1$.
Thus for $\phi(x)\in S^2_\xi(\g_{-1}^*)$ the curve $C_\phi$ 
is linearly equivalent to 
$\ell_\mu+\ell_\nu$, where $(\ell_\mu^2)=(\ell_\nu^2)=-1$,
$(\ell_\mu.\ell_\nu)=1$, so that $(C_\phi^2)=0$.
The Riemann--Roch theorem implies 
$\dim \H^0(X',\O(C_\phi))=2$, hence $C_\phi$ belongs to
a pencil of curves whose generic members are irreducible conics on $X'$.


\medskip

Let us denote by $\g_{-1}^\times$ the open subset of 
$\g_{-1}$ consisting of the 
points with all weight coordinates non-zero. Similarly,
$X'^\times$ denotes the complement to the union of exceptional curves 
in $X'$. 
Since $X={\rm Bl}_M(X')$ is a del Pezzo surface of degree 1,
we have $M\in X'^\times$. Then
$f'^{-1}(M)\subset \g_{-1}^\times$, that is, the coordinates 
of any point in the fibre above $M$ are non-zero.
Let $x_0\in \T'$ be a $k$-point in the fibre over $M$.
For $y\in\g_{-1}^\times$ we let $\frac{y}{x_0}$ denote the
element of the diagonal torus of $\GL(\g_{-1})$ that sends
$x_0$ to $y$.

\ble \label{choice} 
There exists a non-empty open set 
$\Omega\subset(G'/P')_a$ such that for any $y\in \Omega(\ov k)$, 
any root $\mu$ of $\g'$, any weight $\nu$ of $\g_{-1}$ 
and any quadratic polynomial $s(x)$ of weight $0$ with
respect to $H'$ neither of the forms
$p^{\mu}(x)$, $q^{\nu}(x)-x^{\nu} s(x)$
vanishes identically on $\frac{y}{x_0}\T'$.
\ele
{\em Proof} (cf. \cite[Prop. 6.2]{I}, the first statement)
For contradiction assume that
$p^{\mu}(xy/x_0)$ vanishes at every point $(x,y)$ of $\T'\times(G'/P')_a$.
Up to proportionality $p^\mu(y)$ is a unique 
quadratic polynomial in $I((G'/P')_a)$ of weight $\mu$. 
So for any $x\in \T'$ we have
$$p^{\mu}\left(\frac{x}{x_0}y\right)=tp^{\mu}(y)$$
for some $t\in k^*$. Write
$$p^{\mu}(y)=\sum_{\mu_1+\mu_2=\mu}c_{\mu_1,\mu_2}y^{\mu_1}y^{\mu_2}.$$
By symmetry $c_{\mu_1,\mu_2}\neq 0$ whenever $\mu_1+\mu_2=\mu$.
We can choose a point
$x\in\T'(\ov k)$ such that $f'(x)$ belongs to exactly 
one exceptional curve of $X'$. If this curve corresponds to the weight 
$\mu_1$, then $x^{\mu_1}=0$ and $x^{\nu}\not=0$
for any $\nu\not=\mu_1$. It follows that $t=0$, a contradiction.

Now assume that for any $x\in\T'$ we have
$$q^{\nu}\left(\frac{x}{x_0}y\right)-\frac{x^{\nu}}{x_0^{\nu}}y^{\nu}
s\left(\frac{x}{x_0}y\right)\in I((G'/P')_a).$$ 
We choose a point $x\in\T'(\ov k)$ such that $f'(x)$ lies on 
the exceptional curve corresponding to $\nu$ and no other 
exceptional curve of $X'$.
Then $x^{\nu}=0$ is the only vanishing coordinate of $x$.
Since $q^{\nu}(x)$ is not divisible by $x^\nu$ we obtain 
a contradiction with Lemma~\ref{aux1}. QED

\medskip

Let us fix an open set $\Omega$ as in
Lemma \ref{choice}, and pick up a $k$-point $y_0$ in $\Omega^\times$. Define
$$\tilde\T'=\frac{y_0}{x_0}\T', \quad\quad \tilde X'=\tilde\T'/T',\quad\quad
\tilde p(x)=p\left(\frac{y_0}{x_0}x\right).$$
Let $\tilde M$ be the point $f'(y_0)\in \tilde X'$. An obvious
isomorphism $X'\tilde\lra \tilde X'$ sends $M$ to $\tilde M$, so that
$X$ is isomorphic to the blowing-up of $\tilde M$ in $\tilde X'$.

\ble \label{lemma0}
If $\mu\in \hat T'$ is a weight of $\g'$, $\mu\not=\chi_0$,
then the closed subset of $\T'$ given by $\tilde p^{\mu}(x)=0$
is $f'^{-1}(P_\mu)$, where $P_\mu\subset X'$ is
the unique geometrically integral conic passing through $M$ such that 
$[P_\mu]=\mu$.
\ele
{\em Proof} (cf. \cite[Cor. 6.3]{I}) 
We have seen above that $P_\mu$ is a conic such that $[P_\mu]=\mu$.
Now $y_0\in (G'/P')_a$ implies $\tilde p(x_0)=0$, so that $M\in P_\mu$.
The conic $P_\mu$ cannot be reducible since $M$ lies in $X'^\times$. QED

\bco \label{inte}
The orbit $T'y_0$ is the scheme-theoretic intersection
$\tilde\T'\cap (G'/P')_a$. This implies the following relation
among the tangent spaces at $y_0$:
\begin{equation}
\TT _{y_0,(G'/P')_a}\, \cap \, \TT _{y_0,\tilde\T'}\,=\, \TT _{y_0,T'y_0}. 
\label{tan}
\end{equation}
\eco
{\em Proof} (cf. \cite[Cor. 6.4]{I})
We can easily find two weights $\mu$ and $\nu$ such that
the intersection index of the conics $P_\mu$ and $P_\nu$ is 1, that
is, $M$ is the point of intersection of $P_\mu$ and $P_\nu$ with
multiplicity 1. Hence the orbit $T'y_0$ is the scheme-theoretic intersection
of $\tilde\T'$ and the subvariety of $\g_{-1}$
given by $p^\mu(x)=p^\nu(x)=0$. This implies our statement. QED

\bpr \label{sx}
There exists a quadratic form 
$s(x)\in S^2_{\chi_0}(\g_{-1}^*)$ such that
\begin{equation}\label{condition}
s(y_0)=0, \quad \langle y_0,a\rangle +4s(y_0,a)=0\ \text{for any}\ a\in \TT_{y_0,\tilde\T'}.
\end{equation}
It is unique up to addition of a form from the ideal of $\tilde\T'$.
\epr
{\em Proof} We write $\kappa:\tilde X'\to\P^2_k=
\P\big(\H^0(\tilde X',-K_{\tilde X'})^*\big)$ for
the anti-canonical double covering. By Lemma \ref{PQ}
the induced map
$\kappa^*:\TT^*_{\kappa(\tilde M),\P^2}\to\TT^*_{\tilde M,\tilde X'}$,
is an isomorphism.
Since $f':\tilde\T'\to \tilde X'$ is a torsor under $T'$ we have 
$\TT_{\tilde M,\tilde X'}=\TT_{y_0,\tilde\T'}/\TT_{y_0,T'y_0}$,
so the induced map $f'^*:\TT^*_{\tilde M,\tilde X'}\to\TT^*_{y_0,\tilde\T'}$
is identified with the canonical injection
$(\TT_{y_0,\tilde\T'}/\TT_{y_0,T'y_0})^*\to \TT^*_{y_0,\tilde\T'}$.
The morphisms $f'$ and $\kappa$ thus induce the following maps:
$$\TT^*_{\kappa(\tilde M),\P^2} \tilde\lra \TT^*_{\tilde M,\tilde X'}
\tilde\lra (\TT_{y_0,\tilde\T'}/\TT_{y_0,T'y_0})^*\hookrightarrow \TT^*_{y_0,\tilde\T'}.$$
By (\ref{form}) we have
$$\H^0(\tilde X',-K_{\tilde X'})=
S^2_{\chi_0}(\g^*_{-1})/\big(I(\tilde\T')\cap S^2_{\chi_0}(\g^*_{-1})\big).$$
There is a canonical isomorphism
$$\TT^*_{\kappa(\tilde M),\P^2}=\{s\in S^2_{\chi_0}(\g^*_{-1})/
\big(I(\tilde\T')\cap S^2_{\chi_0}(\g^*_{-1})\big)\ 
\text{such that}\ s(y_0)=0\}.$$
Consider the linear form $L\in \TT^*_{y_0,\tilde\T'}$ defined by
$L(a)=\langle y_0,a\rangle$, where $a\in \TT_{y_0,\tilde\T'}$. 
For any $y\in (G'/P')_a$ and any $a\in \TT _{y,(G'/P')_a}$ we have  
$\langle y,a\rangle=0$ by Lemma \ref{sym}.
In particular, $\TT _{y_0,T'y_0}\subset \Ker(L)$, hence $L$ belongs
to the subspace $(\TT_{y_0,\tilde T'}/\TT_{y_0,T'y_0})^*$. 
It is straightforward to check that the map
$f'^*\kappa^*:\TT^*_{\kappa(\tilde M),\P^2}\to \TT^*_{y_0,\tilde\T'}$
sends $s$ to the linear form $s(y_0,a)$, where $a\in \TT^*_{y_0,\tilde\T'}$.
Therefore, there exists a quadratic form 
$s\in S^2_{\chi_0}(\g_{-1}^*)$ satisfying 
(\ref{condition}). Its uniqueness modulo the ideal of $\tilde\T'$
is clear. QED

\medskip

Let us now define
$$\tilde q(x)=q\left(\frac{y_0}{x_0}x\right)-\frac{y_0}{x_0}x\,s\left(\frac{y_0}{x_0}x\right).$$

\ble \label{lemma1}
If $\mu\in \hat T'$ is a weight of $\g_{-1}$, 
then the closed subset of $\T'$ 
given by $\tilde q^{\mu}(x)=0$ is $f'^{-1}(Q_\mu)$,
where $Q_\mu$ is the unique rational curve
with a double point at $M$ and no other singularities, such that
$[Q_\mu]=\chi_0+\mu=-K_{X'}+[\ell_\mu]$. 
\ele
{\em Proof} (cf. \cite[Prop. 6.2]{I}, the second statement)
To check that $M\in Q_\mu$ set $x=x_0$. We have $s(y_0)=0$.
Now $y_0\in(G'/P')_a$ implies $p(y_0)=0$ by Lemma \ref{L1}, and so
$q^\mu(y_0)=0$.

Formula (\ref{q}) and condition (\ref{condition})
show that the derivatives of $\tilde q(x)$ vanish on $\TT _{y_0,\tilde\T'}$.

If the curve $Q_\mu$ is not geometrically integral, the condition
$(Q_\mu.-K_{X'})=3$ implies that $Q_\mu$ is either 
the union of three exceptional curves, or the union of an exceptional curve
and a conic. But $M$ is singular on $Q_\mu$, so $M$ must belong to an
exceptional curve, which is a contradiction. QED

\section{The main theorem} \label{five}

Recall from the introduction that $T\subset\GL(\g)$ is the
extension of the maximal torus $H\subset G$ by the centre of $\GL(\g)$.
The torus $T$ is generated by $T'$ and the 1-dimensional torus
$D\subset\GL(\g)$, whose element $t\in k^*$ acts on $\g_n$ as
multiplication by $t^{n+1}$.
We remind the reader that $X={\rm Bl}_M(X')$ is the del Pezzo surface 
of degree 1 obtained by blowing up the point $M$ on $X'$.
Under the canonical isomorphism $X'\tilde\lra\tilde X'$, the point 
$\tilde M$ in $\tilde X'$ corresponds to $M$ in $X'$.
By the main theorem of \cite{I} we have a universal torsor 
$f':\tilde \T'\to \tilde X'$, where
$\tilde\T'$ is a locally closed subset of $\g_{-1}$.

Let us apply Theorem \ref{B2} to $Z=\tilde\T'$ and the map
$s:S^2(\g_{-1})\to \g_{-2}$ given by $s(x)w$,
where $s(x)$ is the quadratic form as in Proposition \ref{sx}.
In this case $Z_0=T'y_0=f'^{-1}(\tilde M)$ is one $T'$-orbit.
Define $\T=\sZ$. This is a locally closed subset
of $(G/P)_a\subset \g$. By Theorem \ref{B2}
we obtain the following commutative diagram
$$\xymatrix{
\T  \ar[r] \ar[rd]&{\rm Bl}_{T'y_0}(\T') \ar[r] \ar[d]&X \ar[d]^\si\\
&\T' \ar[r]^{f'}& X'}$$
where the horizontal arrows are torsors under tori,
and the vertical arrows are contractions with smooth centres. Since $\exp(x+s(x))v$
is $T'$-equivariant, the torus $T'$ acts on $\T$. The 1-dimensional torus
$D$ acts on $\T$ by construction, hence $T$ acts on $\T$. The fibres
of $f:\T\to X$ are orbits of $T$, hence $\T$, as a composition
of two torsors, is an $X$-torsor under $T$.

Let us recall that $\Pic X$ with the integral bilinear form defined by the
intersection index is identified with the orthogonal direct sum of 
$\Z K_X$ and $Q(\EE_8)=P(\EE_8)$,
where $(K_X)^2=1$, and $Q(\EE_8)$ is equipped with the standard invariant
integral bilinear form multiplied by $-1$, see \cite[Ch. 4]{M}.
If $\beta\in Q(\EE_8)$ is a root of $\g$, we let $\ell_\beta$ be
the exceptional curve on $X$ whose class is $[\ell_\beta]=-K_X+\beta$.
These gives all the
240 exceptional curves on $X$. The intersection index $(\ell_\beta.\ell_\gamma)=
1-(\beta,\gamma)$ for any roots $\beta,\,\gamma\in \EE_8$.

Recall that $\omega\in\EE_8$ is the highest weight of $\g$. 
By Theorem \ref{B2} (iii)
the hyperplane section $y^\omega=0$ of $\T$ is $f^{-1}(\ell_\omega)$,
because $\ell_\omega=\si^{-1}(M)$ is the exceptional divisor of 
$\si:X\to X'$.
By construction, for any root $\beta$ of $\g_{-1}$ the hyperplane section
$y^\beta=0$ of $\T$ is $f^{-1}(\ell_\beta)$. The same is true
if $\beta$ is in $\g'$ or in $\g_1$, by Lemma \ref{lemma0}
and Lemma \ref{lemma1}, respectively.

\medskip

Our next goal is to show that $\T\subset(G/P)_a^{\rm sf}$, where
the latter set consists of stable points for the action of $H$
(i.e. the points whose $H$-orbits in $V$ are closed and have finite stabilisers)
with the additional condition that the stabilisers in $T$ are trivial, 
cf. \cite[Def. 2.5]{I}. By geometric invariant theory there exists
a quasi-projective variety $Y$ and a map 
$(G/P)_a^{\rm sf}\to Y$ which is a torsor under $T$.

If $y\in \g$ denote by $\wt(y)$ the set of roots
$\alpha$ such that $y^{\alpha}\neq 0$ and by $\wt^i(y)$ the set of
roots $\alpha$ of the graded component $\g_i$ such that $y^{\alpha}\neq 0$.
By the Hilbert--Mumford criterion $y$ is stable if and only if
$0$ belongs to the interior 
of the convex hull of $\wt(y)$. The stabiliser of $y$ in $T$ is
trivial if the set $\alpha-\beta$ for all $\alpha,\beta\in \wt(y)$ generates
the root lattice of $\g$.

\ble\label{fat} If $y\in (G/P)_a$ satisfies conditions {\rm (i)} and
{\rm (ii)} below, then $y\in (G/P)_a^{\rm sf}$:

{\rm (i)} if $\mu$ and $\nu$ are roots of $\g'$ and
$(\mu,\nu)=1$, then $\mu\in \wt(y)$ or $\nu\in \wt(y)$;

{\rm (ii)} $\wt^1(y)$ and $\wt^{-1}(y)$ are not empty.
\ele
{\em Proof} First, let us prove that $y$ is stable. We can apply
Prop. 2.4 from \cite{I} to the adjoint representation of $\g'$,
since in the case $\EE_7$ it is a fundamental representation. 
By (i) $\wt^0(y)$ satisfies
the condition of this proposition, and hence $0$ is an interior point
of the convex hull of $\wt^0(y)$. By (ii) the convex hull of
$\wt^0(y)$ is not a face of the convex hull of $\wt(y)$, hence $0$ is
in the interior of the convex hull of $\wt(y)$.

Now let us prove that the stabiliser of $y$ in $T$ is trivial. By the
previous result this stabiliser is finite. By Proposition 2.2 of
\cite{I} the differences $\alpha-\beta$ for all 
$\alpha,\beta\in \wt(y)$ generate
the root lattice of some semisimple Lie subalgebra of $\g$ of rank 8.
By (i) this subalgebra contains $\g'$ and (ii) ensures that it
coincides with $\g$. QED

\ble \label{stab} 
The torsor $\T$ is a Zariski closed subset of $(G/P)_a^{\rm sf}$.
\ele
{\em Proof} First, let us prove that $\T\subset (G/P)_a^{\rm sf}$. 
We use Lemma \ref{fat} and prove that any $y\in\T$
satisfies the conditions (i) and (ii). 
Let $\mu$ and $\nu$ be roots of $\g_0$ such that $(\mu,\nu)=1$.
Then the corresponding exceptional curves
$\ell_\mu$ and $\ell_\nu$ are disjoint since 
$(\ell_\mu.\ell_\nu)=1-(\mu,\nu)=0$. Thus either $\mu\in \wt(y)$
or $\nu\in\wt(y)$, which proves (i).

Assume now that $\wt^1(y)=\emptyset$. Take
any two roots $\mu$ and $\nu$ 
of $\g_1$ such that $(\mu,\nu)=1$. Then as above we have
$\ell_\mu\cap\ell_\nu=\emptyset$, hence either $\mu\in \wt^1(y)$
or $\nu\in\wt^1(y)$, so that $\wt^1(y)$ cannot be empty.
The set $\wt^{-1}(y)$ is non-empty since for any point of $X'$
there exists a exceptional curve on $X'$ that does not contain it.
This proves that $\T\subset (G/P)^{\rm sf}_a$. 

We see that $X$ is a subset of $Y$. Since $X$ is proper, 
$\T=f^{-1}(X)$ is closed in $(G/P)_a^{\rm sf}=f^{-1}(Y)$.
QED

\bthe\label{main} 
For the closed embedding $X\hookrightarrow Y$ constructed above,
$\T=f^{-1}(X)$ is a universal $X$-torsor. Moreover,
the $T$-invariant hyperplane sections of $\T$
defined by the roots of $\g$
are the inverse images of the exceptional curves on $X$.
\ethe
{\em Proof} We know that $\T\to X$ is a torsor under $T$, and we also
know that $(G/P)^{\rm sf}_a\to Y$ is a universal torsor, that is,
its type $\hat T\to \Pic Y$ is an isomorphism.
We pointed out above that if $\beta$ is a root of 
$\g_{-2}\oplus\g_{-1}$,
then $y^\beta=0$ is $f^{-1}(\ell_\beta)$. 
Since $[\ell_\omega]$ and $[\ell_\beta]$ for all roots $\beta$
of $\g_{-1}$ generate the abelian group $\Pic X$, the 
restriction map $\Pic Y\to \Pic X$ is surjective. 
Since the ranks of $\Pic Y$ and $\Pic X$ are the same, 
the restriction map is an isomorphism.
Hence the type $\hat T\to \Pic X$ is an isomorphism.
Moreover, it is easy to see that
this isomorphism sends each root $\beta$ of $\g$
to the class of the corresponding exceptional curve $\ell_\beta$.
The last claim of the theorem is already proved 
for all the roots of $\g$
except the one that spans $\g_2$.
For that root the claim is proved in Lemma \ref{excep} below. QED

\ble\label{excep} 
Let $\T\subset (G/P)_a^{\rm sf}$ be a universal $X$-torsor
whose type $\hat T\tilde\lra \Pic X$
sends each root $\beta$ of $\g$ to the class of the 
corresponding exceptional curve $\ell_\beta\subset X$.
If $\beta$ is a root of $\g$, then the hyperplane section 
$y^\beta=0$ of $\T$ is $f^{-1}(\ell_\beta)$.
\ele
{\em Proof} We first show that $\T$ is not contained
in the hyperplane section $y^\beta=0$.
Let $R$ be the $k$-algebra of regular functions on $(G/P)_a$.
In the proof of Thm. 2.7 of \cite{I} we showed that the codimension
of the complement to $(G/P)_a^{\rm sf}$ in $(G/P)_a$ is at most 2.
Hence $R$ is also the algebra of regular functions on 
$(G/P)_a^{\rm sf}$.
Let $k[\T]$ be the algebra of regular functions on $\T$.
The closed embedding $\T\subset (G/P)_a^{\rm sf}$ gives rise to
a natural surjective homomorphism of $k$-algebras $\Phi: R\to k[\T]$.
The action of $T$ on $R$ and $k[\T]$ equips these algebras with
compatible $\hat T$-gradings:
$$R=\bigoplus_{\chi\in\hat T} R_\chi, 
\quad k[\T]=\bigoplus_{\chi\in\hat T} k[\T]_\chi,$$
where $R_\chi$ (respectively, $k[\T]_\chi$) denotes the $T$-eigenspace
of weight $\chi$. Since $\Phi$ is $T$-equivariant
and surjective we must have $\Phi(R_\chi)=k[\T]_\chi$
for every $\chi\in\hat T$. 
If $\chi=-\ell_\beta$, where $\beta$ is a root of $\g$, 
then $R_{-\ell_\beta}$ is spanned by the weight coordinate $y^\beta$. 
Since $\T$ is a universal $X$-torsor, we have
$k[\T]_{-\ell_\beta}=\H^0(X,\O(\ell_\beta))\cong k$, see (\ref{form}). Thus
$\Phi$ defines an isomorphism of 1-dimensional vector spaces
$R_{-\ell_\beta}\tilde\lra k[\T]_{-\ell_\beta}$, in particular,
$\Phi(y^\beta)\neq 0$, so that the hyperplane section of $\T$
given by $y^\beta=0$ is the inverse image of a curve $C\subset X$.
By assumption, in $\Pic X$ we have $[C]=[\ell_\beta]$, hence 
$C=\ell_\beta$,
because $\ell_\beta$ is the only effective divisor in its class. QED

\noindent Department of Mathematics, University of California,
Berkeley, CA, 94720-3840 USA
\medskip

\noindent serganov@math.berkeley.edu

\bigskip

\noindent Department of Mathematics, South Kensington Campus, 
Imperial College London, SW7 2BZ England, U.K.

\smallskip

\noindent Institute for the Information Transmission Problems, 
Russian Academy of Sciences, 19 Bolshoi Karetnyi, 
Moscow, 127994 Russia
\medskip

\noindent a.skorobogatov@imperial.ac.uk

\end{document}